\documentclass[11pt,letterpaper,reqno]{amsart}
\usepackage[T1]{fontenc}
\usepackage{lmodern}
\usepackage[margin=1.05in]{geometry}
\usepackage{amssymb}
\usepackage[expansion=false]{microtype}
\usepackage[hidelinks]{hyperref}
\newtheorem{theorem}{Theorem}[section]
\newtheorem{lemma}[theorem]{Lemma}
\newtheorem{proposition}[theorem]{Proposition}
\newtheorem{corollary}[theorem]{Corollary}
\theoremstyle{definition}
\newtheorem{definition}[theorem]{Definition}
\newtheorem{example}[theorem]{Example}
\theoremstyle{remark}
\newtheorem{remark}[theorem]{Remark}
\numberwithin{equation}{section}

\newcommand{\C}{\mathbb C}
\newcommand{\D}{\mathbb D}

\newcommand{\id}{\mathrm{id}}
\newcommand{\norm}[1]{\lVert#1\rVert}

\DeclareMathOperator{\Hol}{Hol}
\DeclareMathOperator{\dist}{dist}
\DeclareMathOperator{\diam}{diam}
\DeclareMathOperator{\arctanh}{arctanh}
\DeclareMathOperator{\Log}{Log}
\DeclareMathOperator{\Arg}{Arg}
\DeclareMathOperator{\Rea}{Re}
\DeclareMathOperator{\Ima}{Im}
\newcommand{\doi}[1]{\href{https://doi.org/#1}{\nolinkurl{doi:#1}}}

\title[A Burns--Krantz theorem on convex domains]{A Burns--Krantz theorem for bounded convex domains}
\author[C. Lou]{Cheng Lou}
\author[J. Qiao]{Jianyong Qiao}
\author[H. Wang]{Hongyu Wang}
\author[Y. Zhong]{Yumin Zhong}
\address[All authors]{School of Science, Beijing University of Posts and Telecommunications, Beijing 100876, China, Key Laboratory of Mathematics and Information Networks (Beijing University of Posts and Telecommunications), Ministry of Education, China}
\email{chenglou@bupt.edu.cn;qjy@bupt.edu.cn;wanghyu@bupt.edu.cn;ymzhong@bupt.edu.cn}
\date{}

\begin{document}

\begin{abstract}
We prove two boundary rigidity theorems for holomorphic self-maps of bounded convex domains. Let $\Omega\subset\C^n$, $n\geq1$, and suppose that $\partial\Omega$ is $C^{1,\alpha}$ near $p\in\partial\Omega$ for some $0<\alpha\leq1$. If a holomorphic self-map $f$ satisfies $f(z)=z+o(|z-p|^3)$ as $z\to p$ in one fixed non-tangential approach region at $p$, then $f=\id_\Omega$. 

We also prove a rigidity theorem for bounded convex domains $\Omega$ without assuming boundary smoothness. If a holomorphic self-map $f$ of $\Omega$ satisfies $f(z)-z=O(|z-p|^m)$ as $z\to p\in\partial\Omega$ for some $m>1+2\pi/\theta_p$, where $\theta_p$ is the supremum of sector angles in the open tangent cone at $p$, taken over complex lines through its vertex, then $f=\id_\Omega$, which
proves a conjecture of Zimmer.
\end{abstract}

\maketitle

\section{Introduction}

By Cartan's uniqueness theorem, a holomorphic self-map of a bounded domain is the identity if it fixes an interior point and its differential at that point is the identity. Burns and Krantz~\cite{BurnsKrantz} proved a boundary analogue for smoothly bounded strongly pseudoconvex domains. Let $\Hol(U,V)$ denote the set of holomorphic maps from $U$ to $V$. If $p\in\partial\Omega$ and $f\in\Hol(\Omega,\Omega)$ satisfies
\[
f(z)=z+o(|z-p|^3)\qquad(\Omega\ni z\to p),
\]
then $f=\mathrm{id}_\Omega$. Examples in the unit disc show that rigidity can fail if the exponent $3$ is replaced by any smaller exponent, or if $o(|z-p|^3)$ is replaced by $O(|z-p|^3)$.

Huang~\cite{Huang} then studied boundary rigidity on smoothly bounded convex domains of finite type. If $\ell(p)$ denotes the line type at $p$, his result gives $f=\id_\Omega$ whenever $f(z)=z+o(|z-p|^m)$ as $\Omega\ni z\to p$ for some $m>5\ell(p)$. Huang suggested that the exponent $3$ should suffice in this setting as well~\cite[Remark~3.6]{Huang}.

Zimmer~\cite{Zimmer2022} later removed the finite-type assumption, proving that $f(z)=z+o(|z-p|^4)$ implies $f=\id_\Omega$ on any bounded convex domain with $C^2$ boundary. He further showed that if the boundary is $C^\infty$ and $p$ has finite line type, then the exponent $4$ can be lowered to $4-1/\ell(p)$. In the general $C^2$ convex setting, Zimmer asked whether the exponent $3$ suffices~\cite[Remark~1.6(3)]{Zimmer2022}.

In our first result, it suffices to impose the condition $f(z)-z=o(|z-p|^3)$ in a single non-tangential approach region at $p$.

\begin{theorem}\label{thm:smooth}
Let $\Omega\subset\C^n$ ($n\geq1$) be a bounded convex domain, and let $p\in\partial\Omega$. Suppose that $\partial\Omega$ is of class $C^{1,\alpha}$ near $p$ for some $0<\alpha\leq1$. For $a>1$, define
\[
\Gamma_a(p)=\{z\in\Omega:|z-p|<a\delta_\Omega(z)\},\qquad
\delta_\Omega(z)=\dist(z,\partial\Omega).
\]
If $f\in\Hol(\Omega,\Omega)$ satisfies
\[
f(z)=z+o(|z-p|^3)\qquad(\Gamma_a(p)\ni z\to p)
\]
for some $a>1$, then $f=\id_\Omega$.
\end{theorem}

\begin{remark}
\begin{enumerate}
\item The theorem shows that the exponent $3$ suffices for rigidity on smoothly bounded convex domains of finite type, as Huang suggested in~\cite[Remark~3.6]{Huang}. It also answers Zimmer's question for bounded convex domains with $C^2$ boundary~\cite[Remark~1.6(3)]{Zimmer2022}. In Section~4, we prove the same conclusion under the weaker boundary regularity assumption that $p$ is a Dini-smooth boundary point.
\item Example~\ref{ex:sharp-cubic} gives a nonidentity self-map with $O(|z-p|^3)$ contact, and hence with $o(|z-p|^m)$ contact for every $m<3$. For $\theta=\pi$, Example~\ref{ex:critical-cone} gives a nonidentity self-map with cubic little-$o$ contact at a $C^1$ boundary point.
\end{enumerate}
\end{remark}

Non-tangential analogues of the Burns--Krantz theorem have also been obtained in other geometric settings. Baracco, Zaitsev and Zampieri~\cite{BaraccoZaitsevZampieri} treated the case of germs at strongly pseudoconvex boundary points. Rong~\cite[Theorem~3.2 and Section~4]{RongNontangential} considered strongly linearly convex $C^3$ domains. Further results for special classes of domains can be found in~\cite{FornaessRong,NgRong,Zwonek,KumarPonnusamy,RongCorners}; see also~\cite{RongSurvey}.

Our second boundary result requires no regularity assumption on $\partial\Omega$. For $p\in\partial\Omega$, let $C_p$ be the interior of the real tangent cone of $\Omega$ at $p$; equivalently,
\[
C_p=\bigcup_{s>0}s(\Omega-p)=\bigcup_{s>0}\{s(z-p):z\in\Omega\}.
\]
Define
\begin{equation}\label{eq:opening}
\theta_p=\sup\left\{\theta\in(0,\pi):
\begin{array}{l}
\text{there exists }v\in\C^n,\ |v|=1,\text{ such that }e^{it}v\in C_p\\
\text{for all }|t|\leq\theta/2
\end{array}
\right\}.
\end{equation}
Thus $\theta_p$ measures the largest angular aperture of a sector in $C_p$ along some complex line through the vertex.

Zimmer~\cite[Conjecture~1.8]{Zimmer2022} conjectured that for every boundary point $p$ of a bounded convex domain $\Omega$, there exists an exponent $m_p$, depending only on the tangent cone at $p$, such that any holomorphic self-map $f$ of $\Omega$ satisfying
\[
f(z)-z=o(|z-p|^{m_p})\qquad(\Omega\ni z\to p)
\]
is the identity. The following theorem proves the conjecture and gives an explicit range of sufficient exponents in terms of $\theta_p$.

\begin{theorem}\label{thm:cone}
Let $\Omega\subset\C^n$, $n\geq1$, be a bounded convex domain, and let $p\in\partial\Omega$. If $f\in\Hol(\Omega,\Omega)$ satisfies
\[
f(z)-z=O(|z-p|^m)\qquad(\Omega\ni z\to p)
\]
for some real $m>1+2\pi/\theta_p$, then $f=\id_\Omega$.
\end{theorem}

\begin{remark}
The condition $m>1+2\pi/\theta_p$ is sharp for each fixed value of $\theta_p$. Example~\ref{ex:critical-cone} shows that rigidity can fail at the critical exponent $1+2\pi/\theta_p$ even under a little-$o$ condition.
\end{remark}

The main estimate used to prove Theorem~\ref{thm:smooth} and Theorem~\ref{thm:cone} is the following quantitative form of Cartan's uniqueness theorem. Let $K_\Omega$ and $\kappa_\Omega$ denote the Kobayashi distance and infinitesimal metric, respectively, and let $\rho(A)$ denote the spectral radius of a linear map~$A$.

\begin{theorem}\label{thm:cartan}
For every $n\geq1$ there is a constant $C_n>0$ such that, for every bounded convex domain $\Omega\subset\C^n$, every $f\in\Hol(\Omega,\Omega)$, and all $q,z\in\Omega$,
\[
|f(z)-z|\leq C_n\diam(\Omega)e^{4K_\Omega(q,z)}
\bigl(\kappa_\Omega(q;f(q)-q)+\rho(df_q-I)\bigr).
\]
Here $I$ is the identity map on $\C^n$.
\end{theorem}

\begin{remark}
The exponent $4$ in the exponential factor is optimal; see Example~\ref{ex:sharp-cubic}.
\end{remark}

Christodoulou and Short~\cite{ChristodoulouShort} first proved an inequality of this kind for holomorphic self-maps of the unit disc. Zimmer used it in his proof of boundary rigidity for bounded convex domains with $C^2$ boundary~\cite[Section~5]{Zimmer2022}. Their results motivated both Theorem~\ref{thm:cartan} above and Proposition~\ref{prop:criterion} below.

For $f\in\Hol(\Omega,\Omega)$ and $q\in\Omega$, put $G=f-\id_\Omega$ and set
\[
\varepsilon_G(q)=\kappa_\Omega(q;G(q))+\rho(dG_q),\qquad
D(f;q)=\sup_{\substack{z\in\Omega\\ K_\Omega(q,z)<1}}K_\Omega(z,f(z)).
\]

\begin{proposition}\label{prop:comparison}
For each $n\geq1$, there exist constants $\tau_n>0$ and $L_n\geq1$ such that, for every bounded convex domain $\Omega\subset\C^n$, every $f\in\Hol(\Omega,\Omega)$, and every $q\in\Omega$, one has
\[
L_n^{-1}\varepsilon_G(q)\leq D(f;q)\leq L_n\varepsilon_G(q)
\]
whenever $\min\{\varepsilon_G(q),D(f;q)\}\leq\tau_n$, where $G=f-\id_\Omega$.
\end{proposition}

Since $D(f;q)$ is invariant under biholomorphic changes of domain, Proposition~\ref{prop:comparison} gives a version of the sequence criterion in terms of Kobayashi distance. It holds on any complex manifold biholomorphic to a bounded convex domain; see Corollary~\ref{cor:distance-rigidity}.

\section{Preliminaries}

\noindent\textit{Notation.} Throughout the paper, $\Omega$ is a bounded convex domain in $\C^n$, and $\D$ is the unit disc in $\C$.
\begin{enumerate}
\renewcommand{\labelenumi}{\arabic{enumi}.}
\item We use $|\cdot|$ for the absolute value of a scalar and the Euclidean norm of a vector. The symbols $\operatorname{dist}$ and $\operatorname{diam}$ denote Euclidean distance and diameter, respectively, and $B(q,r)$ denotes the open Euclidean ball with center $q$ and radius $r$.
\item For $q\in\Omega$ and $v\in\C^n\setminus\{0\}$, write
\[
\delta_\Omega(q;v)=\operatorname{dist}\bigl(q,\partial\Omega\cap(q+\C v)\bigr).
\]
\item For a complex-linear map $A:\C^n\to\C^n$, let $\norm A=\sup_{|v|=1}|Av|$ denote its Euclidean operator norm.
\end{enumerate}

\noindent\textit{The Kobayashi metric.} For $z\in\Omega$ and $v\in\C^n$, the infinitesimal Kobayashi metric $\kappa_\Omega(z;v)$ is the infimum of all $a>0$ such that some holomorphic disc $\phi:\D\to\Omega$ satisfies $\phi(0)=z$ and $a\phi'(0)=v$. The associated Kobayashi distance is given by
\[
K_\Omega(z,w)=\inf_\gamma\int_0^1\kappa_\Omega(\gamma(t);\gamma'(t))\,dt,
\]
where $\gamma$ ranges over all piecewise $C^1$ curves in $\Omega$ with $\gamma(0)=z$ and $\gamma(1)=w$.

On the unit disc, we have
\[
\kappa_\D(\zeta;v)=\frac{|v|}{1-|\zeta|^2},
\qquad
K_\D(0,r)=\frac12\log\frac{1+r}{1-r}
\qquad(0\leq r<1).
\]
Both the Kobayashi distance and the infinitesimal Kobayashi metric are non-increasing under holomorphic maps. If $F:\Omega\to\Omega'$ is holomorphic, then
\begin{equation}\label{666}
\kappa_{\Omega'}(F(z);dF_z(v))\leq\kappa_\Omega(z;v),
\qquad
K_{\Omega'}(F(z),F(w))\leq K_\Omega(z,w)
\end{equation}
for all $z,w\in\Omega$ and $v\in\C^n$.

For $q\in\Omega$ and $R>0$, we write
\[
B_K(q,R)=\{z\in\Omega:K_\Omega(q,z)<R\},
\qquad
\overline B_K(q,R)=\{z\in\Omega:K_\Omega(q,z)\leq R\}
\]
for the open and closed Kobayashi balls of radius $R$ centered at $q$, respectively. On a bounded convex domain, $K_\Omega$ is complete and proper, so every closed Kobayashi ball is compact in $\Omega$ \cite[Theorem~2.1]{Zimmer2022}.

Graham obtained the following estimate for the infinitesimal Kobayashi metric on convex domains.

\begin{theorem}[Graham \cite{Graham}]\label{thm:graham}
Let $\Omega\subset\C^n$ be a convex domain. For every $q\in\Omega$ and every $v\in\C^n\setminus\{0\}$,
\[
\frac{|v|}{2\delta_\Omega(q;v)}\leq\kappa_\Omega(q;v)\leq\frac{|v|}{\delta_\Omega(q;v)}.
\]
In particular, for every $q\in\Omega$,
\begin{equation}\label{eq:indicatrix}
\{v\in\C^n:\kappa_\Omega(q;v)<\tfrac12\}\subset\Omega-q=\{z-q:z\in\Omega\}.
\end{equation}
\end{theorem}

\begin{lemma}\label{lem:convex-kobayashi-norm}
Let $\Omega\subset\C^n$ be a bounded convex domain and let $q\in\Omega$. Then $\kappa_\Omega(q;\cdot)$ is a complex norm on $\C^n$.
\end{lemma}

\begin{proof}
Positive definiteness follows from Theorem~\ref{thm:graham}, and complex homogeneity follows from the definition. To prove the triangle inequality, choose holomorphic discs $\phi,\psi:\D\to\Omega$ centered at $q$ with $a\phi'(0)=v$ and $b\psi'(0)=w$, where $a,b>0$. By convexity, $(a\phi+b\psi)/(a+b)$ is a holomorphic disc in $\Omega$ centered at $q$ with derivative $(v+w)/(a+b)$ at $0$. Thus $\kappa_\Omega(q;v+w)\leq a+b$. Taking infima over the admissible $a$ and $b$ gives the triangle inequality.
\end{proof}

\begin{definition}\label{def:complex-geodesic}
A holomorphic map $\phi:\D\to\Omega$ is a complex geodesic if
\[
K_\Omega(\phi(\zeta),\phi(\eta))=K_\D(\zeta,\eta)
\qquad(\zeta,\eta\in\D).
\]
A holomorphic left inverse of $\phi$ is a map $\pi\in\Hol(\Omega,\D)$ satisfying $\pi\circ\phi=\id_\D$.
\end{definition}

For the following form of Lempert's theorem, see also \cite[Theorem~2.7 and Corollary~2.8]{Zimmer2022}.

\begin{lemma}[Lempert \cite{Lempert1981,Lempert1982}]\label{lem:lempert}
Let $\Omega\subset\C^n$ be a bounded convex domain, and let $q,x\in\Omega$ be distinct.
\begin{enumerate}
\renewcommand{\labelenumi}{(\alph{enumi})}
\item If $\Omega$ is strongly convex with $C^\infty$ boundary, there is a unique complex geodesic image containing $q$ and $x$. Every complex geodesic parametrizing this image admits a holomorphic left inverse $\pi:\Omega\to\D$.
\item There exist a complex geodesic $\phi:\D\to\Omega$ and a holomorphic left inverse $\pi:\Omega\to\D$ of $\phi$ such that
\[
\phi(0)=q,\qquad\phi\bigl(\tanh K_\Omega(q,x)\bigr)=x.
\]
\end{enumerate}
In both parts, the left inverse $\pi$ can be chosen so that
\[
\pi^{-1}(\zeta)=\Omega\cap H_\zeta\qquad(\zeta\in\D),
\]
where each $H_\zeta$ is a complex affine hyperplane.
\end{lemma}

The following estimate follows from the Schwarz--Pick lemma (see also \cite[p.~635, (1)]{Harris}).

\begin{lemma}\label{lem:disc-estimate}
Let $g\in\Hol(\D,\D)$ satisfy $g(0)=0$. Then
\[
|g(r)-r|\leq\frac{r(1+r)}{1-r}|g'(0)-1|
\qquad(0\leq r<1).
\]
In particular, if $|g'(0)-1|\leq\varepsilon$, then
\[
|g(r)-r|\leq\frac{2\varepsilon}{1-r}
\qquad(0\leq r<1).
\]
\end{lemma}

\begin{proof}
Let $g(\zeta)=\zeta h(\zeta)$, where $h$ is holomorphic on $\D$
and $h(0)=g'(0)$. Schwarz's lemma gives $|h|\le1$ on $\D$.
If $h$ is nonconstant, then $h(\D)\subset\D$ by the maximum
modulus principle, and the Schwarz--Pick lemma gives
\[
|h(r)-h(0)|\le r|1-\overline{h(0)}h(r)|.
\]
The same inequality holds when $h$ is constant. Using $|h(0)|\le1$,
we obtain
\begin{align*}
|1-h(r)|
\le |1-h(0)|+r|1-\overline{h(0)}h(r)|
\le (1+r)|1-h(0)|+r|1-h(r)|.
\end{align*}
Moving the last term to the left and using
$|g(r)-r|=r|1-h(r)|$ gives
\[
|g(r)-r|\le \frac{r(1+r)}{1-r}|1-g'(0)|.
\]
The second estimate follows from $r(1+r)\le2$.
\end{proof}

\section{A quantitative form of Cartan's uniqueness theorem}

In this section we prove Theorem~\ref{thm:cartan}. We first establish several preliminary results. For a complex norm $N$ on $\C^n$, denote
\[
\norm A_N=\sup_{v\ne0}N(Av)/N(v)
\]
for the induced operator norm.

\begin{lemma}\label{lem:matrix}
Let $N$ be a complex norm on $\C^n$, and let $A$ be a
complex $n\times n$ matrix. Suppose that
\begin{equation}\label{111}
\norm{e^{tA}}_N\leq Ke^{\lambda t}\qquad(t\geq0)
\end{equation}
for some $K\geq1$ and $\lambda\geq0$. Then
\[
\norm A_N\leq\gamma_{n,K}\bigl(\lambda+\rho(A)\bigr),
\]
where $\gamma_{n,K}$ depends only on $n$ and $K$.
\end{lemma}

\begin{proof}
Take $r>\lambda+\rho(A)$ and put $B=A/r$, we have
$\rho(B)<1$.  By \eqref{111}
\begin{equation}\label{222}
\norm{e^{tB}}_N
=\norm{e^{(t/r)A}}_N
\leq Ke^{(\lambda/r)t}
\leq Ke^t.
\end{equation}
Since $\rho(B)<1$, the matrix $2I-B$ is invertible.
We have
\[
\frac{d}{dt}\bigl(e^{-2t}e^{tB}\bigr)
=-(2I-B)e^{-2t}e^{tB}.
\]
Integrating over $[0,T]$ therefore gives
\[
(2I-B)\int_0^T e^{-2t}e^{tB}\,dt
=I-e^{-2T}e^{TB}.
\]
Multiplying both sides on the left by
$R=(2I-B)^{-1}$, we obtain
\[
\int_0^T e^{-2t}e^{tB}\,dt
=R-Re^{-2T}e^{TB}.
\]
By \eqref{222}, we have
\[
\norm{e^{-2T}e^{TB}}_N
\leq Ke^{-T}\longrightarrow0
\qquad(T\to\infty).
\]
Letting $T\to\infty$, we conclude that
\[
R=\int_0^\infty e^{-2t}e^{tB}\,dt.
\]
Furthermore, \eqref{222} yields
\[
\begin{aligned}
\norm R_N
\leq\int_0^\infty
\norm{e^{-2t}e^{tB}}_N\,dt
\leq K\int_0^\infty e^{-t}\,dt
=K.
\end{aligned}
\]
In particular,
\[
\norm R_N
\leq\int_0^\infty\norm{e^{-2t}e^{tB}}_N\,dt
\leq K\int_0^\infty e^{-t}\,dt
=K.
\]

It remains to estimate the inverse of $R$. Let
$\mu_1,\ldots,\mu_n$ denote the eigenvalues of $B$,
counted with algebraic multiplicity. The eigenvalues of $R$
are $\eta_i=(2-\mu_i)^{-1}$. The condition $|\mu_i|<1$
implies
$
\frac13<|\eta_i|<1.
$
Write the characteristic polynomial of $R$ as
\[
p_R(s)=\prod_{i=1}^n(s-\eta_i)
=\sum_{j=0}^n a_js^j.
\]
The bounds on the eigenvalues give
\[
|a_0|\geq 3^{-n},
\qquad
|a_j|\leq\binom nj,\quad 1\leq j\leq n.
\]
Indeed, $a_0$ is, up to sign, the product of all the
$\eta_i$, whereas each $a_j$ is a sum of
$\binom nj$ products of $n-j$ eigenvalues.
By the Cayley--Hamilton theorem,
\[
a_0I+a_1R+\cdots+a_nR^n=0.
\]
As $R$ is invertible, this identity can be rewritten as
\[
R^{-1}=-\frac1{a_0}\sum_{j=1}^n a_jR^{j-1}.
\]
Using the preceding coefficient estimates and $\norm R_N\leq K$,
we arrive at
\[
\begin{aligned}
\norm{R^{-1}}_N
\leq 3^n\sum_{j=1}^n\binom nj K^{j-1}
\leq 3^n(2^n-1)K^{n-1}.
\end{aligned}
\]

Finally, $B=2I-R^{-1}$ and $A=rB$, so
\[
\norm A_N
\leq r\bigl(2+3^n(2^n-1)K^{n-1}\bigr).
\]
Since this holds for every $r>\lambda+\rho(A)$,
letting $r\rightarrow\lambda+\rho(A)$ completes
the proof.
\end{proof}
For $q\in\Omega$, let $\norm{\cdot}_q$ denote the operator norm induced by $\kappa_\Omega(q;\cdot)$. Thus, for a complex-linear map $A:\C^n\to\C^n$,
\[
\norm A_q=\sup_{v\ne0}\frac{\kappa_\Omega(q;Av)}{\kappa_\Omega(q;v)}.
\]

\begin{lemma}\label{lem:fixed-correction}
Let $f\in\Hol(\Omega,\Omega)$ and put $G=f-\id_\Omega$. For any $q\in\Omega$, the map
\begin{equation}\label{eq:fixed-correction}
Y(z)=G(z)-G(q)-4\kappa_\Omega(q;G(q))(z-q)
\end{equation}
can be written as $(1+4\kappa_\Omega(q;G(q)))(h-\id_\Omega)$ for some $h\in\Hol(\Omega,\Omega)$ with $h(q)=q$. Moreover,
\[
\norm{dY_q}_q\leq\norm{dG_q}_q+4\kappa_\Omega(q;G(q)),\qquad
|G(z)-Y(z)|\leq 4\kappa_\Omega(q;G(q))\operatorname{diam}(\Omega)
\qquad(z\in\Omega).
\]
\end{lemma}

\begin{proof}
Fix $q\in\Omega$, set $a=4\kappa_\Omega(q;G(q))$. Differentiating both sides of the equation \eqref{eq:fixed-correction}, we obtain  $$dY_q=dG_q-aI.$$ Since $\norm I_q=1$, we have $\norm{dY_q}_q\leq\norm{dG_q}_q+a$. If $a=0$, then $f(q)=q$ and $Y=G$, so we may take $h=f$. Suppose now that $a>0$ and set
\[
b_0=q-\frac{G(q)}a,\qquad h(z)=\frac{f(z)+ab_0}{1+a}.
\]
Since $\kappa_\Omega(q;b_0-q)=1/4$, we have $b_0\in\Omega$ by \eqref{eq:indicatrix} and thus $h\in\Hol(\Omega,\Omega)$ by convexity. Also we obtain
\[
(1+a)(h(z)-z)=G(z)-a(z-b_0)=Y(z).
\]
As $Y(q)=0$, this gives $h(q)=q$. Finally it also gives
\[
|G(z)-Y(z)|=a|z-b_0|\leq a\operatorname{diam}(\Omega).
\qedhere
\]
\end{proof}

For $f\in\Hol(\Omega,\Omega)$ and $G=f-\id_\Omega$, we also set
\[
\varepsilon_G^{\mathrm{op}}(q)=\kappa_\Omega(q;G(q))+\norm{dG_q}_q
\qquad(q\in\Omega).
\]

\begin{proposition}\label{prop:spectral}
There is a constant $c_n>0$, depending only on $n$, such that for every $f\in\Hol(\Omega,\Omega)$, with $G=f-\id_\Omega$,
\[
\norm{dG_q}_q\leq c_n\varepsilon_G(q)
\qquad(q\in\Omega).
\]
Consequently,
\[
\varepsilon_G(q)\leq\varepsilon_G^{\mathrm{op}}(q)\leq(c_n+1)\varepsilon_G(q)
\qquad(q\in\Omega).
\]
\end{proposition}

\begin{proof}
Fix $q\in\Omega$ and let $a=4\kappa_\Omega(q;G(q))$, $Y$, and $h$ be as in Lemma~\ref{lem:fixed-correction}. Since $h(q)=q$, we have $\norm{dh_q}_q\leq1$ by \eqref{666}. The identities $Y=(1+a)(h-\id_\Omega)$ and $dY_q=dG_q-aI$ give
\[
dG_q=(1+a)dh_q-I.
\]
Then for $t\geq0$,
\[
\norm{e^{t\,dG_q}}_q=e^{-t}\norm{e^{t(1+a)dh_q}}_q
\leq e^{-t}e^{t(1+a)}=e^{at}.
\]
Applying Lemma~\ref{lem:matrix} with $N=\kappa_\Omega(q;\cdot)$, $K=1$, and $\lambda=a$, we obtain
\[
\norm{dG_q}_q\leq\gamma_{n,1}\bigl(a+\rho(dG_q)\bigr)\leq4\gamma_{n,1}\varepsilon_G(q).
\]
Thus we may take $c_n=4\gamma_{n,1}$. The inequality $\rho(dG_q)\leq\norm{dG_q}_q$ gives $\varepsilon_G(q)\leq\varepsilon_G^{\mathrm{op}}(q)$, and
\[
\varepsilon_G^{\mathrm{op}}(q)=\kappa_\Omega(q;G(q))+\norm{dG_q}_q\leq(c_n+1)\varepsilon_G(q).
\qedhere
\]
\end{proof}

\begin{remark}
The estimate in Proposition~\ref{prop:spectral} can fail for a general holomorphic map. For example, the matrix $A=\left(\begin{smallmatrix}0&1\\0&0\end{smallmatrix}\right)$ has $\rho(A)=0$ and $\norm A=1$. If $q\in\Omega\subset\C^2$, the map $X(z)=A(z-q)$ satisfies
\[
\kappa_\Omega(q;X(q))+\rho(dX_q)=0,\qquad dX_q=A\ne0.
\]
\end{remark}

\begin{remark}
For each $q\in\Omega$, monotonicity of the Kobayashi metric under the inclusions
\[
B(q,\delta_\Omega(q))\subset\Omega\subset B(q,\operatorname{diam}(\Omega))
\]
gives
\[
\norm{dG_q}_q\leq\frac{\operatorname{diam}(\Omega)}{\delta_\Omega(q)}\norm{dG_q}.
\]
This comparison introduces a factor of $\delta_\Omega(q)^{-1}$ when passing from the Euclidean operator norm to the operator norm at $q$. Since $\rho(dG_q)\leq\norm{dG_q}$, Cauchy's estimate bounds the spectral-radius term in $\varepsilon_G(q)$ directly. We use this observation in the proof of Theorem~\ref{thm:smooth}.
\end{remark}

\begin{lemma}\label{lem:resolvent}
Let $f\in\Hol(\Omega,\Omega)$ and put $G=f-\id_\Omega$. For every $t\geq0$ and $w\in\Omega$, there is a unique point $R_t(w)\in\Omega$ satisfying
\begin{equation}\label{eq:resolvent}
R_t(w)-tG(R_t(w))=w.
\end{equation}
The map $R_t:\Omega\to\Omega$ is holomorphic and injective. We call $R_t$ the resolvent of $G$. If $q\in\Omega$ and $G(q)=0$, then $R_t(q)=q$ and
\[
d(R_t)_q=(I-t\,dG_q)^{-1}.
\]
\end{lemma}

\begin{proof}
For $t=0$, we have $R_0=\id_\Omega$. For $t>0$, the equation $y-tG(y)=w$ takes the form
\[
y+t(y-f(y))=w.
\]
By Reich and Shoikhet \cite[Proposition~4.2(1)]{ReichShoikhet}, this equation has a unique solution in $\Omega$ that depends holomorphically on $w$.

Equation \eqref{eq:resolvent} gives $(\id_\Omega-tG)\circ R_t=\id_\Omega$, so $R_t$ is injective. If $G(q)=0$, then $q$ solves the equation with $w=q$. Uniqueness gives $R_t(q)=q$, and differentiation at $q$ yields
\[
(I-t\,dG_q)d(R_t)_q=I.
\qedhere
\]
\end{proof}

\begin{lemma}\label{lem:exit}
Let $f\in\Hol(\Omega,\Omega)$ and put
$G=f-\id_\Omega$. Suppose that $q\in\Omega$
and $G(q)=0$. For every $z\in\Omega$ with
$G(z)\ne0$, there exist $t>0$ and $w\in\Omega$
such that
\[
G(z)=\frac{z-w}{t},
\]
and
\[
K_\Omega(q,w)=K_\Omega(q,z)+1,
\qquad
\frac1t\leq
e^{4K_\Omega(q,z)+2}\norm{dG_q}_q.
\]
\end{lemma}

\begin{proof}
Fix $z\in\Omega$ with $G(z)\ne0$ and set
$k=K_\Omega(q,z)$ and $U=B_K(q,k+1)$.
Consider the ray $z-sG(z)$, $s\geq0$,
let $t>0$ be the first time the ray leaving $U$
and put $w=z-tG(z)$. Since the closure of $U$
is compact in $\Omega$, 
$w\in\partial U\subset\Omega$, and  we have
\[
K_\Omega(q,w)=k+1.
\]
On the other hand, consider the equations
$$z-tG(z)=w.$$ Lemma~\ref{lem:resolvent} gives
$R_t(w)=z$ and $R_t(q)=q$.

Set $r=\tanh K_\Omega(q,w)$.
By Lemma~\ref{lem:lempert}, there exists a
complex geodesic $\phi$ with $\phi(0)=q$
and $\phi(r)=w$, together with a holomorphic
left inverse $\pi$. By \eqref{666},
\[
\kappa_\Omega(q;\phi'(0))=1,
\qquad
|d\pi_q(v)|\leq\kappa_\Omega(q;v).
\]
Then $g=\pi\circ R_t\circ\phi$
is a holomorphic self-map of $\D$ with
$g(0)=0$ and $g(r)=\pi(z)$. By \eqref{666},
\[
K_\D(0,g(r))
=K_\D(\pi(q),\pi(z))\leq k.
\]
In particular, $|g(r)|\leq\tanh k$.

Since $R_t(q)=q$, also by \eqref{666},
we have
\(
\norm{d(R_t)_q}_q\leq1.
\)
Thus Lemma~\ref{lem:resolvent} gives
\[
d(R_t)_q-I=t\,dG_q\,d(R_t)_q
\quad\text{and}\quad
\norm{d(R_t)_q-I}_q
\leq t\norm{dG_q}_q.
\]
Since $\pi\circ\phi=\id_\D$, the chain rule gives
\[
|g'(0)-1|
=\bigl|d\pi_q\bigl((d(R_t)_q-I)\phi'(0)\bigr)\bigr|
\leq t\norm{dG_q}_q.
\]
Together with $|g(r)|\leq\tanh k$,
Lemma~\ref{lem:disc-estimate} yields
\[
r-\tanh k\leq |r-g(r)|
\leq\frac{r(1+r)}{1-r}\,
t\norm{dG_q}_q.
\]
Since
\[
r=\tanh(k+1)
=\frac{e^{2k+2}-1}{e^{2k+2}+1},
\]
we have \(
\frac{1+r}{1-r}=e^{2k+2}
\)and
\[
r-\tanh k
=\frac{2(e^2-1)e^{2k}}
{(e^{2k+2}+1)(e^{2k}+1)}.
\]
Thus
\[
\begin{aligned}
\frac1t
\leq\frac{r(1+r)}
{(1-r)(r-\tanh k)}
\norm{dG_q}_q
=\frac{e^2(e^{2k}+1)(e^{2k+2}-1)}
{2(e^2-1)}\norm{dG_q}_q
\leq e^{4k+2}\norm{dG_q}_q.
\end{aligned}
\]
\end{proof}

\begin{proof}[Proof of Theorem~\ref{thm:cartan}]

Put $G=f-\id_\Omega$ and $a=4\kappa_\Omega(q;G(q))$.
Let $Y$ and $h$ be given by Lemma~\ref{lem:fixed-correction}, so that
$Y=(1+a)(h-\id_\Omega)$ and $h(q)=q$.
If $h(z)\ne z$, Lemma~\ref{lem:exit} gives $s>0$ and $w\in\Omega$
such that $$h(z)-z=(z-w)/s$$ and
$$s^{-1}\leq e^{4K_\Omega(q,z)+2}\norm{dh_q-I}_q.$$ Since
$|z-w|\leq\operatorname{diam}(\Omega)$, we have
\[
|Y(z)|=(1+a)|h(z)-z|
\leq(1+a)\operatorname{diam}(\Omega)e^{4K_\Omega(q,z)+2}
       \norm{dh_q-I}_q
=\operatorname{diam}(\Omega)e^{4K_\Omega(q,z)+2}\norm{dY_q}_q.
\]
By Lemma~\ref{lem:fixed-correction}, we obtain
\begin{align*}
|G(z)|
&\leq |Y(z)|+|G(z)-Y(z)|
\leq \operatorname{diam}(\Omega)e^{4K_\Omega(q,z)+2}
       \bigl(\norm{dG_q}_q+a\bigr)
       +a\operatorname{diam}(\Omega)\\
&\leq \operatorname{diam}(\Omega)e^{4K_\Omega(q,z)+2}
       \bigl(\norm{dG_q}_q+8\kappa_\Omega(q;G(q))\bigr),
\end{align*}
Then by Proposition~\ref{prop:spectral},
$\norm{dG_q}_q\leq c_n\varepsilon_G(q)$. Since
$\kappa_\Omega(q;G(q))\leq\varepsilon_G(q)$, we may take
$C_n=e^2(c_n+8)$.
\end{proof}

\section{Two boundary rigidity theorems}

In this section we prove Theorems~\ref{thm:smooth} and~\ref{thm:cone}.

\begin{proposition}\label{prop:criterion}
Let $f\in\Hol(\Omega,\Omega)$ and put $G=f-\id_\Omega$. If there exist $o\in\Omega$ and a sequence $(q_j)\subset\Omega$ such that
\[
e^{4K_\Omega(o,q_j)}\varepsilon_G(q_j)\longrightarrow 0,
\]
then $f=\id_\Omega$. In particular, if there exist $o\in\Omega$ and a sequence  $\{s_j\} \in\mathbb{R}^{+}$ satisfying $s_j\leq \delta_\Omega(q_j)$ such that
\[
s_j^{-1}e^{4K_\Omega(o,q_j)}\sup_{B(q_j,s_j)}|G|\longrightarrow 0.
\]then $f=\id_\Omega$.
\end{proposition}

\begin{proof}
	Fix $z\in\Omega$. By Theorem~\ref{thm:cartan} and the
	triangle inequality for the Kobayashi distance, we have
	\[
	\begin{aligned}
		|G(z)|
		\leq C_n\diam(\Omega)e^{4K_\Omega(q_j,z)}
		\varepsilon_G(q_j)
		\leq C_n\diam(\Omega)e^{4K_\Omega(o,z)}
		e^{4K_\Omega(o,q_j)}\varepsilon_G(q_j).
	\end{aligned}
	\]
	The right-hand side tends to zero as $j\to\infty$,
	and hence $G(z)=0$.
	
	Suppose that $B(q,s)\subset\Omega$, by \eqref{666},
	\[
	\kappa_\Omega(q;G(q))
	\leq\kappa_{B(q,s)}(q;G(q))
	=\frac{|G(q)|}{s}.
	\]
	On the other hand, for every Euclidean unit vector $v$,
	Cauchy's estimate applied to the holomorphic map
	$\zeta\mapsto G(q+\zeta v)$ gives
	\[
	|dG_qv|\leq\frac1s\sup_{B(q,s)}|G|.
	\]
	Note that $\rho(dG_q)\leq\norm{dG_q}$, taking the supremum over all such $v$,
we obtain
	\[
	\begin{aligned}
		\varepsilon_G(q)
		=\kappa_\Omega(q;G(q))+\rho(dG_q)
		\leq\frac{|G(q)|}{s}+\norm{dG_q}
		\leq\frac{2}{s}\sup_{B(q,s)}|G|.
	\end{aligned}
	\]
	which completes the proof.
\end{proof}

We will prove Theorem~\ref{thm:smooth} under a boundary regularity assumption slightly weaker than $C^{1,\alpha}$ smoothness. We call $p\in\partial\Omega$ a Dini-smooth boundary point if $\partial\Omega$ is $C^1$ near $p$ and its unit normal is Dini continuous there. In particular, every $C^{1,\alpha}$ boundary point with $0<\alpha\leq1$ is Dini smooth. The following lemma is a consequence of Nikolov and Thomas~\cite{NikolovThomas2025} and the local estimate of Nikolov and Andreev~\cite{NikolovAndreev}.

\begin{lemma}\label{lem:normal-distance}
Let $p$ be a Dini-smooth boundary point of $\Omega$ and let $\nu_p$ be the outward unit normal at $p$. Set $p_r=p-r\nu_p$. Then $p_r\in\Omega$ for all sufficiently small $r>0$, and
\[
\frac{\delta_\Omega(p_r)}{r}\rightarrow 1
\qquad (r\rightarrow0).
\]
For every $o\in\Omega$, there is a constant $C_o$ such that
\[
K_\Omega(o,p_r)\leq C_o+\tfrac12\log(1/r)
\]
for all sufficiently small $r>0$.
\end{lemma}
\begin{proof}
	For the proofs, we refer to \cite[Proposition~7]{NikolovThomas2025} and  \cite[Theorem~7]{NikolovAndreev}. In fact, the first assertion requires only that $\partial\Omega$ be $C^1$-smooth.
\end{proof}

\begin{proof}[Proof of Theorem~\ref{thm:smooth}]
Put $G=f-\id_\Omega$ and take $0<\eta<(a-1)/(a+1)$ for $a>1$.
By Lemma~\ref{lem:normal-distance}, $\delta_\Omega(p_r)/r\to1$,
so $B(p_r,\eta r)\subset\Omega$ for sufficiently small $r>0$.
For $z\in B(p_r,\eta r)$, we have $|z-p|\leq(1+\eta)r$ and
$\delta_\Omega(z)\geq\delta_\Omega(p_r)-\eta r$. Hence
\[
\frac{|z-p|}{\delta_\Omega(z)}
\leq\frac{1+\eta}{\delta_\Omega(p_r)/r-\eta}<a.
\]
 Thus $B(p_r,\eta r)\subset\Gamma_a(p)$,
and the assumption gives
\[
\sup_{B(p_r,\eta r)}|G|=o(r^3).
\]
Fix $o\in\Omega$, by Lemma~\ref{lem:normal-distance} we have
$e^{4K_\Omega(o,p_r)}=O(r^{-2})$ and 
\[
\frac{e^{4K_\Omega(o,p_r)}}{\eta r}
\sup_{B(p_r,\eta r)}|G|\longrightarrow0.
\]
Thus we have
$f=\id_\Omega$ by Proposition~\ref{prop:criterion}
\end{proof}

In the proof of Theorem~\ref{thm:cone}, we place the centers of
the Euclidean balls on the image of the real radius of an analytic
disc approaching $p$.

\begin{lemma}\label{lem:sector-disc}
	For every $0<\theta<\theta_p$, there exist a unit vector
	$v\in\C^n$ and  $a>0$ such that the map $H:\D\rightarrow \Omega$ given by
	\[
	H(\zeta)=p+a(1-\zeta)^{\theta/\pi}v
	\]
is holomorphic, here
	\[
	(1-\zeta)^{\theta/\pi}
	=\exp\left(\frac{\theta}{\pi}\Log(1-\zeta)\right),
	\]
	where $\Log$ denotes the principal branch of the logarithm on the
	right half-plane.
\end{lemma}

\begin{proof}
Since $\Omega$ is convex and $p\in\overline\Omega$, we have
$t(\Omega-p)\subset\Omega-p$ for $0<t\leq1$. Thus the open sets
$r^{-1}(\Omega-p)$ increase to $C_p$ as $r\downarrow0$.

By the definition of $\theta_p$, there is a unit vector $v$ such
that the arc $\{e^{it}v:|t|\leq\theta/2\}$ lies in $C_p$.
The arc is compact, so it is contained in $R^{-1}(\Omega-p)$
for some $R>0$. Scaling toward $p$ now gives
\[
p+\lambda v\in\Omega
\qquad(0<|\lambda|<R,\ |\operatorname{Arg}\lambda|<\theta/2).
\]
For $\zeta\in\D$, we have $\Rea(1-\zeta)>0$ and
$|1-\zeta|<2$. Thus, 
\[
\left|\arg(1-\zeta)^{\theta/\pi}\right|<\frac{\theta}{2}
\quad \text{and} \quad
\left|(1-\zeta)^{\theta/\pi}\right|<2^{\theta/\pi}.
\]
Thus we can choose $a>0$ sufficiently small so that
$a2^{\theta/\pi}<R$ which shows that
$H(\D)\subset\Omega$.

\end{proof}
\begin{proof}[Proof of Theorem~\ref{thm:cone}]
Put $G=f-\id_\Omega$. Choose $\theta$ with
$2\pi/(m-1)<\theta<\theta_p$ and set $\beta=\theta/\pi$.
Let $H:\D\to\Omega$ be the holomorphic map given by
Lemma~\ref{lem:sector-disc}.
For $0\leq t<1$, take
\[
q_t=H(t)\quad \text{and}\quad s_t=(1-t)^\beta\delta_\Omega(q_0).
\]
Since $q_t-p=(1-t)^\beta(q_0-p)$, by convexity we have
\[
B(q_t,s_t)
=p+(1-t)^\beta\bigl(B(q_0,\delta_\Omega(q_0))-p\bigr)
\subset\Omega.
\]
Then for $z\in B(q_t,s_t)$, we have
\[
|z-p|\leq |q_t-p|+s_t
=(1-t)^\beta\bigl(|q_0-p|+\delta_\Omega(q_0)\bigr).
\]
Thus by the assumption we can deduce that
$$\sup_{B(q_t,s_t)}|G|=O((1-t)^{\beta m})\qquad(t\rightarrow1).$$ 
On the other hand, by \eqref{666}, we have
\[
e^{4K_\Omega(q_0,q_t)}
\leq e^{4K_\D(0,t)}
=\left(\frac{1+t}{1-t}\right)^2.
\]
It follows that
\[
\frac{e^{4K_\Omega(q_0,q_t)}}{s_t}\sup_{B(q_t,s_t)}|G|
=O\bigl((1-t)^{\beta(m-1)-2}\bigr)\rightarrow0,
\]
since $\beta(m-1)>2$. By Proposition~\ref{prop:criterion} we completes the proof.
\end{proof}

\section{Rigidity in terms of the Kobayashi distance}

In this section, we prove Proposition~\ref{prop:comparison} and use it to obtain the biholomorphically invariant form of Proposition~\ref{prop:criterion} stated in Corollary~\ref{cor:distance-rigidity}.

\begin{lemma}\label{lem:metric-comparison}
	For all $q,z\in\Omega$ and $v\in\C^n$,
	\[
	\kappa_\Omega(q;z-q)\leq e^{2K_\Omega(q,z)}-1,
	\qquad
	\kappa_\Omega(z;v)\leq
	2e^{2K_\Omega(q,z)}\kappa_\Omega(q;v).
	\]
	If $x\in\Omega$ satisfies
	$\kappa_\Omega(q;x-q)\leq1/4$, then
	\[
	\tfrac13\kappa_\Omega(q;v)
	\leq\kappa_\Omega(x;v)
	\leq3\kappa_\Omega(q;v)
	\qquad(v\in\C^n).
	\]
\end{lemma}

\begin{proof}
	The first two estimates are immediate when $q=z$.
	Suppose that $q\ne z$. By Theorem~\ref{thm:graham}
	and \cite[Remark~2.7]{ZimmerGromov}, we have
	\[
	\kappa_\Omega(q;z-q)
	\leq\frac{|z-q|}{\delta_\Omega(q;z-q)}
	\leq e^{2K_\Omega(q,z)}-1.
	\]
	
	For the second estimate, let $a$ be the boundary point
	on the ray from $q$ through $z$, and put
	\[
	\lambda=\frac{|z-a|}{|q-a|}.
	\]
	The affine map $A(w)=a+\lambda(w-a)$ maps $\Omega$
	into itself by convexity. Since $A(q)=z$ and
	$dA_q(v)=\lambda v$, the contraction property of
	the Kobayashi metric gives
	\[
	\lambda\kappa_\Omega(z;v)
	=\kappa_\Omega(A(q);dA_q(v))
	\leq\kappa_\Omega(q;v).
	\]
	By \cite[Lemma~2.6]{ZimmerGromov},
	$\lambda^{-1}\leq e^{2K_\Omega(q,z)}$. Consequently,
	\[
	\kappa_\Omega(z;v)
	\leq e^{2K_\Omega(q,z)}\kappa_\Omega(q;v).
	\]

	Now suppose that
	$\kappa_\Omega(q;x-q)\leq1/4$.
	By~\eqref{eq:indicatrix}, the map
	\[
	\zeta\longmapsto q+2\zeta(x-q)
	\]
	takes $\D$ into $\Omega$. Since it sends $0$ to $q$
	and $1/2$ to $x$, the distance-decreasing property
	gives
	\[
	K_\Omega(q,x)\leq K_\D(0,1/2)
	=\arctanh(1/2)=\tfrac12\log3.
	\]
	Applying the preceding estimate in both directions,
	we obtain
	\[
	\begin{aligned}
		\kappa_\Omega(x;v)
		&\leq e^{2K_\Omega(q,x)}\kappa_\Omega(q;v)
		\leq3\kappa_\Omega(q;v),\\
		\kappa_\Omega(q;v)
		&\leq e^{2K_\Omega(q,x)}\kappa_\Omega(x;v)
		\leq3\kappa_\Omega(x;v).
	\end{aligned}
	\]
	These inequalities imply the desired comparison.
	\qedhere
\end{proof}

\begin{corollary}\label{cor:small-distance}
If $q,x\in\Omega$ and $\kappa_\Omega(q;x-q)<1/2$, then
\[
K_\Omega(q,x)\leq\arctanh\bigl(2\kappa_\Omega(q;x-q)\bigr).
\]
In particular, if $\kappa_\Omega(q;x-q)\leq1/4$, then $K_\Omega(q,x)\leq4\kappa_\Omega(q;x-q)$.
\end{corollary}

\begin{proof}
The assertion is immediate when $x=q$. Otherwise, put $b=\kappa_\Omega(q;x-q)$. The disc
\[
\zeta\longmapsto q+\frac{\zeta(x-q)}{2b}
\]
takes values in $\Omega$ by~\eqref{eq:indicatrix}, since the Kobayashi norm at $q$ of its displacement from $q$ is $|\zeta|/2<1/2$. The disc takes $0$ to $q$ and $2b$ to $x$, so $K_\Omega(q,x)\leq\arctanh(2b)$ by contraction. For $b\leq1/4$, the inequality $\arctanh(2b)\leq4b$ gives the second assertion.
\end{proof}

\begin{lemma}\label{lem:jet-from-displacement}
There is an absolute constant $B\geq1$ such that, for $f\in\Hol(\Omega,\Omega)$, $G=f-\id_\Omega$, and $q\in\Omega$,
\[
D(f;q)\leq1\implies\varepsilon_G^{\mathrm{op}}(q)\leq BD(f;q).
\]
\end{lemma}

\begin{proof}
Suppose $D(f;q)\leq1$. If $\kappa_\Omega(q;x-q)\leq1/4$ then by~\eqref{eq:indicatrix}, $x\in\Omega$. By Corollary~\ref{cor:small-distance} we have  $K_\Omega(q,x)\leq\arctanh(1/2)<1$, so $K_\Omega(x,f(x))\leq D(f;q)$. By Lemma~\ref{lem:metric-comparison},
\[
\kappa_\Omega(q;G(x))\leq3\kappa_\Omega(x;G(x))
\leq3\bigl(e^{2D(f;q)}-1\bigr)\leq3(e^2-1)D(f;q).
\]

Fix $v$ with $\kappa_\Omega(q;v)=1$, the map $\zeta\mapsto G(q+\zeta v)$ is holomorphic for $|\zeta|<1/2$ by~\eqref{eq:indicatrix}. By Cauchy's estimate on \(|\zeta|=1/4\), we have
\[
\kappa_\Omega(q;dG_qv)
\leq 4\sup_{|\zeta|=1/4}\kappa_\Omega(q;G(q+\zeta v))
\leq 12(e^2-1)D(f;q).
\]
Taking the supremum over $\kappa_\Omega(q;v)=1$ and we obtain
$\norm{dG_q}_q\leq 12(e^2-1)D(f;q)$.

Since $K_\Omega(q,f(q))\leq D(f;q)\leq1$, applying Lemma~\ref{lem:metric-comparison}
and we have
\[
\begin{aligned}
	\kappa_\Omega(q;G(q))
	\leq e^{2K_\Omega(q,f(q))}-1
	\leq e^{2D(f;q)}-1
	\leq (e^2-1)D(f;q).
\end{aligned}
\]
Combining these two estimates, we may take
$B=13(e^2-1)$.
\end{proof}

To obtain the converse estimate for $D(f;q)$, we use the
resolvent identity to bound $\kappa_\Omega(z;G(z))$.

\begin{lemma}\label{lem:intrinsic-value}
	Let $f\in\Hol(\Omega,\Omega)$ and $G=f-\id_\Omega$. For all $q,z\in\Omega$,
	\[
	\kappa_\Omega(z;G(z))\leq (4e^4+6)e^{8K_\Omega(q,z)}\varepsilon_G^{\mathrm{op}}(q).
	\]
	If $G(q)=0$, then
	\[
	\kappa_\Omega(z;G(z))\leq e^{8K_\Omega(q,z)+4}\norm{dG_q}_q.
	\]
\end{lemma}

\begin{proof}

Denote $k=K_\Omega(q,z)$. Suppose first that $G(q)=0$.
The desired estimate is immediate if $G(z)=0$.
Otherwise, Lemma~\ref{lem:exit} gives $t>0$ and
$w\in\Omega$ such that $G(z)=\frac{z-w}{t}$ and 
\[
\frac1t\leq e^{4k+2}\norm{dG_q}_q,\qquad
K_\Omega(q,w)=k+1.
\]

Then we have
	\[
	K_\Omega(z,w)\leq K_\Omega(q,z)+K_\Omega(q,w)=2k+1.
	\]
	By Lemma~\ref{lem:metric-comparison},
	\[
	\kappa_\Omega(z;G(z))=\frac{\kappa_\Omega(z;w-z)}{t}
	\leq\frac{e^{2K_\Omega(z,w)}-1}{t}\leq e^{8k+4}\norm{dG_q}_q.
	\]

	For a general $G$, put
	$a=4\kappa_\Omega(q;G(q))$ and let
	$Y=(1+a)(h-\id_\Omega)$ be as in
	Lemma~\ref{lem:fixed-correction}. Since $h(q)=q$,
	we have
	\[
	\kappa_\Omega(z;Y(z))
	\leq e^{8k+4}\norm{dY_q}_q.
	\]
	By Lemma~\ref{lem:fixed-correction},
	\[
	\norm{dY_q}_q
	\leq\norm{dG_q}_q+a
	\leq4\varepsilon_G^{\mathrm{op}}(q).
	\]
	
	On the other hand, since
	$G(z)-Y(z)=G(q)+a(z-q)$, the triangle inequality and 
	Lemma~\ref{lem:metric-comparison} shows that
	\[
	\begin{aligned}
		\kappa_\Omega(z;G(q))+a\kappa_\Omega(z;z-q)
		\leq2e^{2k}\kappa_\Omega(q;G(q))
		+a(e^{2k}-1)
		\leq6e^{2k}\varepsilon_G^{\mathrm{op}}(q).
	\end{aligned}
	\]
	
	The triangle inequality then gives
	\[
	\kappa_\Omega(z;G(z))
	\leq\kappa_\Omega(z;Y(z))+\kappa_\Omega(z;G(z)-Y(z))
	\leq(4e^{8k+4}+6e^{2k})\varepsilon_G^{\mathrm{op}}(q)
	\leq (4e^4+6)e^{8k}\varepsilon_G^{\mathrm{op}}(q). \qedhere
	\]
\end{proof}

\begin{proof}[Proof of Proposition~\ref{prop:comparison}]
	Put $G=f-\id_\Omega$. Let $c_n$ be the constant in
	Proposition~\ref{prop:spectral}, and let $B$ be the
	constant in Lemma~\ref{lem:jet-from-displacement}.
	Set
	\[
	L_n=\max\{B,4(4e^4+6)e^8(c_n+1)\}.
	\]
	For every $z\in B_K(q,1)$, Lemma~\ref{lem:intrinsic-value}
	and Proposition~\ref{prop:spectral} yield
	\[
	4\kappa_\Omega(z;G(z))
	\leq 4(4e^4+6)e^8\varepsilon_G^{\mathrm{op}}(q)
	\leq L_n\varepsilon_G(q).
	\]
	If $\varepsilon_G(q)\leq L_n^{-1}$, then
	$\kappa_\Omega(z;G(z))\leq1/4$ throughout $B_K(q,1)$.
	By Corollary~\ref{cor:small-distance}, we obtain
	\[
	K_\Omega(z,f(z))
	\leq4\kappa_\Omega(z;G(z))
	\leq L_n\varepsilon_G(q).
	\]
	Taking the supremum over $z\in B_K(q,1)$ and we have
	$
	D(f;q)\leq L_n\varepsilon_G(q).
	$
	Conversely, if $D(f;q)\leq1$, Lemma~\ref{lem:jet-from-displacement}
	gives
	\[
	\varepsilon_G(q)\leq B D(f;q)
	\leq L_nD(f;q).
	\]

	Now set $\tau_n=L_n^{-2}$. If
	$\varepsilon_G(q)\leq\tau_n$, then
	\[
	D(f;q)\leq L_n\varepsilon_G(q)
	\leq L_n^{-1}\leq1.
	\]
	Hence the estimate
	$
	\varepsilon_G(q)\leq L_nD(f;q)
	$
	also applies. Conversely, if $D(f;q)\leq\tau_n$, then
	\[
	\varepsilon_G(q)\leq L_nD(f;q)
	\leq L_n^{-1},
	\]
	and consequently
	$
	D(f;q)\leq L_n\varepsilon_G(q).
$
	Thus both inequalities hold whenever either
	$\varepsilon_G(q)$ or $D(f;q)$ is at most $\tau_n$.
	
\end{proof}

For a Kobayashi hyperbolic complex manifold $M$, $f\in \Hol(M,M)$ and $q\in M$, define
\[
D_M(f;q)=\sup_{\substack{z\in M\\ K_M(q,z)<1}}K_M(z,f(z)).
\]

\begin{corollary}\label{cor:distance-rigidity}
	Let $M$ be biholomorphic to a bounded convex domain in $\C^n$ and  $f\in\Hol(M,M)$. If there exist a point $o\in M$ and a sequence $\{q_j\}\subset M$ such that
	\[
	e^{4K_M(o,q_j)}D_M(f;q_j)\longrightarrow0,
	\]
	then $f=\id_M$.
\end{corollary}

\begin{proof}
	Let $\Phi:M\to\Omega$ be a biholomorphism onto a bounded
	convex domain $\Omega\subset\C^n$, and define
	\[
	\widetilde f=\Phi\circ f\circ\Phi^{-1},
	\qquad
	\widetilde o=\Phi(o),
	\qquad
	\widetilde q_j=\Phi(q_j).
	\]
	Since biholomorphisms preserve the Kobayashi distance,
	we have
	\[
	K_\Omega(\widetilde o,\widetilde q_j)
	=K_M(o,q_j),
	\qquad
	D(\widetilde f;\widetilde q_j)=D_M(f;q_j).
	\]
	Thus we have
	\[
	e^{4K_\Omega(\widetilde o,\widetilde q_j)}
	D(\widetilde f;\widetilde q_j)\longrightarrow0.
	\]

By Proposition~\ref{prop:comparison} we have
\[
\begin{aligned}
	e^{4K_\Omega(\widetilde o,\widetilde q_j)}
	\varepsilon_{\widetilde G}(\widetilde q_j)
	\leq L_n e^{4K_\Omega(\widetilde o,\widetilde q_j)}
	D(\widetilde f;\widetilde q_j)
	\rightarrow0.
\end{aligned}
\]
It follows from Proposition~\ref{prop:criterion} that
$\widetilde f=\id_\Omega$
and hence $f=\id_M$.

\end{proof}

\section{Examples}

In this section we start from the same cubic perturbation of the identity
on the unit disc. The first example uses it to show that the
contact condition in Theorem~\ref{thm:smooth} and the exponent $4$ in
Theorem~\ref{thm:cartan} are sharp. In the second example,
we conjugate the disc map by a conformal map to show the rigidity can fail under little-$o$ contact at the critical order for each $0<\theta_p\leq\pi$. When $\theta=\pi$,
this construction also shows that the boundary regularity assumption in Theorem~\ref{thm:smooth} cannot be weakened to $C^1$.

\begin{example}\label{ex:sharp-cubic}
	The polynomial $g_0(\zeta)=\zeta-\frac{1}{10}(\zeta-1)^3$ is a nonidentity holomorphic self-map of $\D$ (see \cite{BurnsKrantz,KrantzSurvey}). Let $\Omega=\D^n$ and $F(z_1,z')=(g_0(z_1),z')$. When $n=1$, we omit the second coordinate block throughout this example. The boundary of $\Omega$ is given by $|z_1|=1$ and is smooth near $p=(1,0)$. For $z\in\Omega$,
	\begin{align*}
		|F(z)-z|=\frac{1}{10}|z_1-1|^3\leq\frac{1}{10}|z-p|^3.   
	\end{align*}
	Thus rigidity can fail under the little-$o$ condition of any order less than three, or even under a cubic big-$O$ condition. In fact, for $q_r=(r,0)$ with $0<r<1$, we have
	\[
	\delta_\Omega(q_r)=|q_r-p|=1-r,\qquad
	|F(q_r)-q_r|=\frac{1}{10}|q_r-p|^3.
	\]
	Since $q_r\in\Gamma_a(p)$ for every $a>1$, the cubic little-$o$ condition fails.
	
	It also shows that the exponent $4$ in Theorem~\ref{thm:cartan} is sharp. Let $G=F-\id_\Omega$. The holomorphic map $\zeta\mapsto(\zeta,0)$ and the first-coordinate projection are both holomorphic, and their composition on $\D$ is $\id_\D$. Contraction in both directions gives
	\[
	K_\Omega(0,q_r)=K_\D(0,r)=\frac12\log\frac{1+r}{1-r},\qquad
	\kappa_\Omega(q_r;(v,0))=\frac{|v|}{1-r^2}\quad(v\in\C).
	\]
	Since $G(q_r)=(\frac{1}{10}(1-r)^3,0)$ and
	$\rho(dG_{q_r})=\frac{3}{10}(1-r)^2$, we obtain
	\begin{align*}
		\varepsilon_G(q_r)=\dfrac{(1-r)^2}{10}\left(3+\frac{1}{1+r}\right).
	\end{align*}
	
	For every $\lambda<4$, this gives
	\begin{align*}
		e^{\lambda K_\Omega(0,q_r)}\varepsilon_G(q_r)
		=O\bigl((1-r)^{2-\lambda/2}\bigr)\longrightarrow0
		\qquad(r\rightarrow1).
	\end{align*}
	So an estimate as in Theorem~\ref{thm:cartan} with using the exponent $\lambda$ in place of $4$ will therefore imply $F=id_\Omega$, which leads a contradiction.
\end{example}

\begin{example}\label{ex:critical-cone}
	For every $n\geq1$ and $0<\theta\leq\pi$, there exist a bounded convex domain $\Omega\subset\C^n$, a point $p\in\partial\Omega$ with $\theta_p=\theta$, and a nonidentity holomorphic self-map $F$ of $\Omega$ such that
	\begin{align*}
		F(z)-z=o\bigl(|z-p|^{1+2\pi/\theta}\bigr)
		\qquad(\Omega\ni z\to p).
	\end{align*}
	When $\theta=\pi$, the boundary can be chosen $C^1$ near $p$.
\end{example}

\begin{proof}
	Let $\beta=\theta/\pi\in(0,1]$. For $\zeta\in\D$ we let $w=1-\zeta$ and define
	\begin{align*}
		L(w)=\frac{3}{\beta}-\Log w,\qquad
		\Psi_\beta(\zeta)=w^\beta L(w),\qquad
		\Omega_\beta=\Psi_\beta(\D),
	\end{align*}
	
	using the principal logarithm on the right half-plane. For $\beta=1$, this is a variant of the logarithmic example in \cite[Example~2]{NikolovPflugThomas}.
	
	Firstly we prove convexity, let $S=\beta L(w)-1$. On $\overline{\D}\setminus\{1\}$, we have $|w|\leq2$ and $\Rea w>0$, so
	\begin{align*}
		\Rea S=2-\beta\log|w|\geq2-\log 2>1
	\end{align*}
	\begin{align*}
		\Rea(wS)=(\Rea w)(\Rea S)+\beta(\Ima w)\Arg w>0.
	\end{align*}
	In the second inequality, $(\Ima w)\Arg w\geq0$. Differentiating gives $\Psi_\beta'(\zeta)=-w^{\beta-1}S$ and
	\begin{align*}
		1+\frac{\zeta\Psi_\beta''(\zeta)}{\Psi_\beta'(\zeta)}
		=\beta\left(1-\frac{1}{S}\right)+\frac{1-\beta}{w}+\frac{\beta}{wS}.
	\end{align*}
	The first term has positive real part since $\Rea S>1$, and the other two also have nonnegative real part.
	The convexity criterion therefore shows that $\Psi_\beta$ is univalent and $\Omega_\beta$ is convex. Since $w^\beta L(w)\to0$ as $w\to0$, the map extends continuously to $\overline{\D}$ with $\Psi_\beta(1)=0$, so $\Omega_\beta$ is bounded. Since $w^\beta L(w)\ne0$ in $\D$, we have $0\in\partial\Omega_\beta$.
	
	Let $C_0$ be the open tangent cone of $\Omega_\beta$ at origin, $w=\rho e^{i\varphi}$ and $u=\frac{3}{\beta}-\log\rho$. Therefore, we see that $|\varphi|<\pi/2$ and $u\beta>1$ by the definition of $w$. Note that $\Psi_\beta$ satisfies that $$ \operatorname{Arg}\Psi_\beta(\zeta) = \beta\varphi-\arctan\left(\frac{\varphi}{u}\right), $$ where the arctangent has the sign of $\varphi$ and its absolute value is at most $\frac{|\varphi|}{u}\le\beta|\varphi|$. Hence $|\operatorname{Arg}\Psi_\beta(\zeta)| \le\beta|\varphi|<\frac{\beta\pi}{2}$. It follows that $\Omega_\beta$, and hence $C_0$, is contained in the sector $$ S_\beta = \{z\in\mathbb{C}\setminus\{0\}: |\operatorname{Arg}z|<\beta\pi/2\}. $$ For fixed $\varphi$ , when $\rho>0$ is sufficiently small, the point $1-\rho e^{i\varphi}$ lies in $\mathbb{D}$ and its image has argument tending to $\beta\varphi$ as \(\rho\) tends to zero through positive values. Therefore, the arguments of points in $\Omega_\beta$ come arbitrarily close to both endpoints of $(-\beta\pi/2,\beta\pi/2)$. Since $\Omega_\beta$ is connected and lies in $S_\beta$, its set of arguments is the whole interval. As a result we have $C_0=S_\beta$.
	
	Since $\Omega_\beta$ is bounded and convex, 
	its boundary is a Jordan curve. 
	Carath\'eodory's theorem gives that the continuous extension of $\Psi_\beta$ is a homeomorphism from $\overline{\D}$ onto $\overline{\Omega_\beta}$.
	
	For $\beta=1$, let $\gamma(t)=\Psi_1(e^{it})$ for the boundary parametrization near $0$. 
	For $t\neq 0$ near $0$, we have $\gamma'(t)=-ie^{it}(2-\Log(1-e^{it}))$ and then $|\gamma'(t)|=O(\log(\dfrac{1}{|t|}))$, which is integrable near $0$. The unit tangent satisfies
	\begin{align*}
		\dfrac{\gamma'(t)}{|\gamma'(t)|}
		=\dfrac{ie^{it}\Psi_1'(e^{it})}{|ie^{it}\Psi_1'(e^{it})|}
		=-ie^{it}\dfrac{2-\Log(1-e^{it})}{|2-\Log(1-e^{it})|}
		\rightarrow -i
		\qquad(t\to0^\pm).
	\end{align*}
	The last equation follows because the real part of $2-\Log(1-e^{it})$ tends to infinity and the imaginary part of it is bounded.
	
	Now we use signed arclength $s$ to parametrize this arc, with $s=0$ at $\gamma(0)$. For $s\ne0$, the derivative of this parametrization is the unit tangent. The difference quotient defining the derivative at $s=0$ can be expressed as the average of that tangent over the interval $[0,s]$. Consequently, its limit as $s\rightarrow 0$ also exists and equals $-i$. This gives a regular $C^1$ boundary parametrization near $0$.
	
	Let $g_0$ be as in Example~\ref{ex:sharp-cubic} and let $F_\beta=\Psi_\beta\circ g_0\circ\Psi_\beta^{-1}$. This is a nonidentity holomorphic self-map of $\Omega_\beta$. For $z=\Psi_\beta(\zeta)$, we have $|z|=|w|^\beta|L(w)|$. Since $\Rea L(w)>1/\beta$, it follows that $w\to0$ as $z\to0$.
	
	Each point $\eta$ on the segment from $\zeta$ to $g_0(\zeta)$ can be written as $\eta=\zeta+sw^3/10$ with $0\leq s\leq1$. 
	Hence we have
	\[
	1-\eta=w\left(1-\frac{s}{10}w^2\right).
	\]
	In particular, $\dfrac{|1-\eta|}{|w|}=1+O(|w|^2)$ uniformly on the segment. Also we know the segment from $w$ to $1-\eta$ is contained in the right half-plane, 
	has length at most $|w|^3/10$ and stays at distance at least $9|w|/10$ from $0$ for small $|w|$. By integrating the derivative of $L$ along this segment we have $L(1-\eta)=L(w)+O(|w|^2)$. Thus we have, uniformly for $\eta$ on this segment,
	\[
	|\Psi_\beta'(\eta)|=O\bigl(|w|^{\beta-1}|L(w)|\bigr).
	\]
	Integrating along the segment from $\zeta$ to $g_0(\zeta)$ yields
	\[
	|F_\beta(z)-z|=O \bigl(|w|^{\beta+2}|L(w)| \bigr).
	\]
	Since $|z|^{1+2/\beta}=|w|^{\beta+2}|L(w)|^{1+2/\beta}$, we have
	\[
	\frac{|F_\beta(z)-z|}{|z|^{1+2/\beta}}=O\bigl(|L(w)|^{-2/\beta}\bigr)\longrightarrow 0.
	\]
	This gives the counterexample in dimension one.
	
	For $n>1$, let $\Omega=\Omega_\beta\times\D^{n-1}$, $p=(0,0)$ and $F(z_1,z')=(F_\beta(z_1),z')$. The open tangent cone at $p$ is $S_\beta\times\C^{n-1}$. Each vector in this cone has a nonzero first component, whose rotations can remain in $S_\beta$ only over an interval of length at most $\beta\pi$. The direction $(1,0)$ admits every smaller opening, so $\theta_p=\beta\pi=\theta$. The desired contact estimate follows from $|z_1|\leq|z-p|$. When $\beta=1$, the boundary near $p$ is the product of the planar $C^1$ boundary near $0$ with a neighborhood of $0$ in $\C^{n-1}$. It is therefore a $C^1$ hypersurface near $p$.
\end{proof}

\end{document}